\documentclass[pre,aps]{revtex4}
\usepackage{graphicx,epsf}
\usepackage{amssymb}
\usepackage{amsmath}
\usepackage{latexsym}
\usepackage[T1]{fontenc}
\usepackage[francais]{babel}

\def\be{\begin{equation}}
\def\ee{\end{equation}}
\def\bq{\begin{eqnarray}}
\def\eq{\end{eqnarray}}
\def\beq{\begin{eqnarray*}}
\def\eeq{\end{eqnarray*}}

\begin{document}

\title{Sur la d\'{e}termination du p\'{e}rim\`{e}tre de l'ovale \`{a} huit centres}

\author{Jean-Marc Ginoux$^{1}$, Jean-Claude Golvin$^2$}

\affiliation{$^1$Laboratoire LIS, CNRS, UMR 7020, Universit\'{e} de Toulon, BP 20132, F-83957 La Garde cedex, France}

\affiliation{$^2$Directeur de Recherches au CNRS, Universit\'{e} de Bordeaux III Michel de Montaigne, France}

\begin{abstract}
Plusieurs recherches ont permis de montrer qu'un ovale \`{a} huit centres co\"{i}ncide presque parfaitement avec l'ellipse construite sur les m\^{e}mes axes et peut \^{e}tre consid\'{e}r\'{e} comme une repr\'{e}sentation de celle-ci \`{a} condition de choisir convenablement les rayons des arcs de cercles qui le composent. Le calcul de son p\'{e}rim\`{e}tre se r\'{e}duit alors \`{a} la simple somme d'arcs de cercles. Pourtant, il ne nous semble pas que ce calcul, qui pourrait s'av\'{e}rer utile, n'ait jamais encore \'{e}t\'{e} effectu\'{e} ni publi\'{e}. L'objet de cette note est donc de pr\'{e}senter une d\'{e}monstration g\'{e}om\'{e}trique de la d\'{e}termination du p\'{e}rim\`{e}tre de l'ovale \`{a} huit centres.\\
{\it Pour citer cet article~: J.-M. Ginoux, J.-C. Golvin, C. R. Acad. Sci. Paris, Ser. I 356 (2018) 1195-1202.}
\end{abstract}

\maketitle

\section{Introduction}
\label{}

Dans son c\'{e}l\`{e}bre trait\'{e} sur les Coniques, compos\'{e} de huit livres et dont seulement quatre nous sont parvenus, Appollonius de Perge (ca. 240 av. J.-C.) \'{e}tudia la parabole, l'hyperbole et l'ellipse \cite{Decorps}. Alors que les principales propri\'{e}t\'{e}s de l'ellipse furent d\'{e}couvertes assez rapidement, la d\'{e}termination de son p\'{e}rim\`{e}tre devait rester, pendant plusieurs si\`{e}cles, inaccessible aux g\'{e}om\`{e}tres. Aussi, d\`{e}s l'Antiquit\'{e} Grecque puis Romaine et tout au long des si\`{e}cles qui suivirent de nombreuses approximations du p\'{e}rim\`{e}tre de l'ellipse furent propos\'{e}es par de c\'{e}l\`{e}bres astronomes ou math\'{e}maticiens. En 1609, Johannes Kepler pr\'{e}senta une formule (\ref{eq1}) consistant \`{a} approcher la circonf\'{e}rence de l'ellipse par celle de deux cercles dont les rayons \'{e}taient respectivement \'{e}gaux \`{a} la moyenne g\'{e}om\'{e}trique et \`{a} la moyenne arithm\'{e}tique des longueurs de son demi grand axe $a$ et de son demi petit axe $b$ :

\begin{equation}
L\left(a, b \right) = \pi \left( a + b \right)
\label{eq1}
\end{equation}

On remarque que si l'ellipse se r\'{e}duit \`{a} un cercle de rayon $R$, on a : $a=b=R$ et l'on retrouve la valeur bien connue du p\'{e}rim\`{e}tre du cercle. Au XVIII\textsuperscript{e} si\`{e}cle, avec l'av\`{e}nement du calcul infinit\'{e}simal, le comte Giulio Fagnano, un math\'{e}maticien italien autodidacte, introduisit une s\'{e}rie de transformations analytiques ing\'{e}nieuses jetant les bases de la th\'{e}orie des \textit{int\'{e}grales elliptiques}. Ses travaux furent publi\'{e}s en 1750 sous la forme de deux volumes intitul\'{e}s \textit{Produzioni Matematiche}. L'ann\'{e}e suivante, le c\'{e}l\`{e}bre math\'{e}maticien et physicien suisse L\'{e}onard Euler, auquel on avait demand\'{e} d'examiner les deux trait\'{e}s de Fagnano, g\'{e}n\'{e}ralisa ses r\'{e}sultats et fournit la fameuse formule de l'int\'{e}grale elliptique (\ref{eq2}) permettant de calculer la valeur du p\'{e}rim\`{e}tre d'une ellipse en fonction de son excentricit\'{e} $\varepsilon = c/a$ avec $c = \sqrt{a^2-b^2}$ :

\begin{equation}
L\left(a, b \right) = 4a \int_0^{\frac{\pi}{2}} \sqrt{1 - \varepsilon^2 \sin^2\left( t\right)} dt
\label{eq2}
\end{equation}

Cette int\'{e}grale elliptique (\ref{eq2}), dite de ``deuxi\`{e}me esp\`{e}ce'', ne pouvant s'exprimer qu'au moyen d'un d\'{e}veloppement en s\'{e}rie, de nouvelles approximations furent r\'{e}guli\`{e}rement propos\'{e}es jusqu'\`{a} nos jours. On peut citer pour m\'{e}moire celle du math\'{e}maticien italien Giuseppe Peano \cite{Peano} en 1889 ainsi que celles fournies par le math\'{e}maticien indien Srinivasa Ramanujan au d\'{e}but du XX\textsuperscript{e} si\`{e}cle \cite{Berndt}. Parall\`{e}lement, d\`{e}s l'Antiquit\'{e} Grecque diff\'{e}rents instruments, comme l'ellipsographe de Proclus, furent imagin\'{e}s pour dessiner une ellipse \cite{Raynaud}. De m\^{e}me, plusieurs m\'{e}thodes, comme celle dite du ``jardinier'', furent d\'{e}velopp\'{e}es. Du point de vue technique, la construction de monuments ou d'\'{e}l\'{e}ments architecturaux de forme apparemment elliptique, comme le c\'{e}l\`{e}bre Colis\'{e}e \`{a} Rome, incita les ing\'{e}nieurs \`{a} rechercher un moyen de repr\'{e}senter une ellipse le plus fid\`{e}lement possible. D'apr\`{e}s les arch\'{e}ologues qui ont \'{e}tudi\'{e} depuis ces \'{e}difices, il semble que la m\'{e}thode employ\'{e}e par ces b\^{a}tisseurs consistait \`{a} approcher la repr\'{e}sentation graphique d'une ellipse par une courbe polycentrique ou plus exactement par un ovale \`{a} plusieurs centres. Ainsi, en 1682, Antoine Desgodetz, architecte du roi Louis XIV, \'{e}tudia, \`{a} la demande de Colbert, les plus c\'{e}l\`{e}bres monuments de Rome. Il \'{e}crivit \`{a} propos du Colis\'{e}e \cite{Desgodetz} :

\begin{quote}
``J'ai trouv\'{e} l'ovale qui forme le mur (\ldots) pour cela j'ai s\'{e}par\'{e} la circonf\'{e}rence en plusieurs arcs que j'ai trouv\'{e} pouvoir \^{e}tre trac\'{e}s au compas par huit centres.''
\end{quote}

Dans son analyse concernant la construction des amphith\'{e}\^{a}tres romains, J.-C. Golvin a montr\'{e} que certains d'entre eux sont construits \`{a} partir d'ovales \`{a} quatre centres \cite{Golvin}. Mark Wilson Jones \cite{Wilson} qui consid\`{e}re que le plan du Colis\'{e}e est un ovale \`{a} huit centres a \'{e}crit :

\begin{quote}
``By contrast the oval is a compound form, approximating to the ellipse by means of segmental arcs of different radius which meet where they share the same tangent. The simplest oval has four arcs drawn from two pairs of focal points on the main axes. More arcs and focal points may be added on this principle; the greater their number the closer the oval resembles the ellipse.''
\end{quote}

Ainsi, il appara\^{i}t que le trac\'{e} de l'ellipse peut \^{e}tre approch\'{e} par un ovale \`{a} huit centres ou plus \cite{Herrera, Rosin1999, Rosin2004}. Dans le cas de l'ovale \`{a} huit centres, la m\'{e}thode de construction graphique peut \^{e}tre simplifi\'{e}e en se limitant au trac\'{e} du quart de l'ellipse, les trois autres quarts se d\'{e}duisant par simple sym\'{e}trie. Elle consiste alors \`{a} d\'{e}terminer les centres et les rayons de trois cercles permettant de tracer trois arcs de cercles qui formeront ainsi une approximation circulaire par morceaux du quart de l'ellipse. Jusqu'\`{a} pr\'{e}sent, certains chercheurs, comme C. Trevisan \cite{Trevisan}, consid\'{e}raient que dans le cas de l'ovale \`{a} huit centres :

\begin{quote}
``Non solo, infatti, sussistono infinite combinazioni per il raggio minore ($R_1$) e quello maggiore ($R_3$), come nell'ovale a quattro centri; ma anche altrettante infinite possibilit\`{a} di scelta per il raggio intermedio ($R_2$).''
\end{quote}

\begin{quote}
``Non seulement cela, en effet, il y a des combinaisons infinies pour le rayon mineur ($R_1$) et le majeur ($R_3$), comme dans l'ovale \`{a} quatre centres; mais aussi autant de possibilit\'{e}s infinies de choix pour le rayon interm\'{e}diaire ($R_2$). ''
\end{quote}

Cependant, cette difficult\'{e}, qui r\'{e}side dans la d\'{e}termination des centres et des rayons des trois arcs de cercles, avait apparemment \'{e}t\'{e} lev\'{e}e en 1908 par F. Honey \cite{Honey}. Dans son tr\`{e}s court article, il propose d'utiliser les cercles ($C_1$) et ($C_2$), respectivement de centres $e$ et $g$, de rayons $r$ et $R$, osculateurs aux sommets de coordonn\'{e}es $(a,0)$ et $(0,b)$ de l'ellipse ($E$) pour d\'{e}terminer le centre $k$ du \textit{cercle interm\'{e}diaire} (voir Fig. 1). Pour y parvenir, il d\'{e}finit deux cercles auxiliaires que nous appellerons par la suite cercles \textit{mineur} et \textit{majeur}. Il consid\`{e}re ensuite que l'un des deux points d'intersection de ces deux cercles est le centre $k$ du \textit{cercle interm\'{e}diaire} auquel il attribue le rayon $p$ \'{e}gal \`{a} ``la moyenne arithm\'{e}tique ou \`{a} la demi-somme des demi-axes''. N\'{e}anmoins, si Honey pr\'{e}sente une m\'{e}thode permettant d'obtenir pr\'{e}cis\'{e}ment les centres et les rayons de ces trois cercles et donc de dessiner ``un ovale \`{a} huit centres qui co\"{i}ncide avec l'ellipse construite sur les m\^{e}mes axes'', il ne fournit pas les angles aux centres $\gamma$, $\beta$ et $\delta$ de ces trois rayons qui lui auraient permis de calculer les longueurs de ces trois arcs de cercles. Cette note a donc pour objet de pr\'{e}senter une d\'{e}monstration g\'{e}om\'{e}trique conduisant \`{a} l'expression g\'{e}n\'{e}rale des valeurs de ces trois angles aux centres en fonction des longueurs des demi-axes majeur et mineur de l'ellipse. La connaissance de ces trois angles associ\'{e}e aux valeurs des trois rayons fournis par Honey permet ainsi de d\'{e}terminer le p\'{e}rim\`{e}tre de l'ovale \`{a} huit centres et donc d'approximer la valeur du p\'{e}rim\`{e}tre de l'ellipse. Ce r\'{e}sultat pourrait s'av\'{e}rer utile pour les g\'{e}om\`{e}tres et les arch\'{e}ologues. On montrera notamment que, dans le cas du Colis\'{e}e, le p\'{e}rim\`{e}tre de l'ovale \`{a} huit centres fournit une excellente approximation du p\'{e}rim\`{e}tre de ce monument calcul\'{e} au moyen de l'int\'{e}grale elliptique.

\section{D\'{e}termination des angles aux centres}
\label{}

Dans son article, Honey \cite{Honey} propose d'utiliser pour le \textit{cercle mineur} (en bleu sur la Fig. 1) le centre $e(a - b^2/a,0)$ du petit cercle osculateur ($C_1$) et le rayon :

\begin{equation}
|ek| = \frac{\left(a-b\right)\left( a + 2b \right)}{2a} = p - r
\label{eq3}
\end{equation}

o\`{u} $p$ repr\'{e}sente la moyenne arithm\'{e}tique ou la demi-somme des demi-axes : $p=\dfrac{a+b}{2}$ et $r = \dfrac{b^2}{a}$ le rayon de courbure du petit cercle osculateur ($C_1$).Pour le \textit{cercle mineur} (en rouge sur la Fig. 1), il emploie le centre $g(0, b - a^2/b)$ du grand cercle osculateur ($C_2$) et le rayon

\begin{equation}
|gk| =\frac{\left(a-b\right)\left( b + 2a \right)}{2b} = R - p
\label{eq4}
\end{equation}

o\`{u} $R = a^2/b$ est le rayon de courbure du grand cercle osculateur ($C_2$).\\

Honey \cite{Honey} d\'{e}finit enfin le centre du \textit{cercle interm\'{e}diaire} comme le point $k$ d'intersection des deux cercles \textit{majeur} et \textit{mineur} et choisit pour son rayon la valeur $p$. D\`{e}s lors, il y a deux mani\`{e}res d'aborder ce probl\`{e}me : l'une alg\'{e}brique, l'autre g\'{e}om\'{e}trique. Bien que nous ayons r\'{e}solu la question de ces deux mani\`{e}res diff\'{e}rentes, nous avons choisi de ne pr\'{e}senter que la d\'{e}monstration g\'{e}om\'{e}trique car elle est bas\'{e}e sur une construction r\'{e}alis\'{e}e uniquement \`{a} la r\`{e}gle et au compas. Nous tenons la d\'{e}monstration alg\'{e}brique \`{a} la disposition des lecteurs qui peuvent l'obtenir par une simple demande.

\begin{figure}[htbp]
\centerline{\includegraphics[width=15cm,height=15cm]{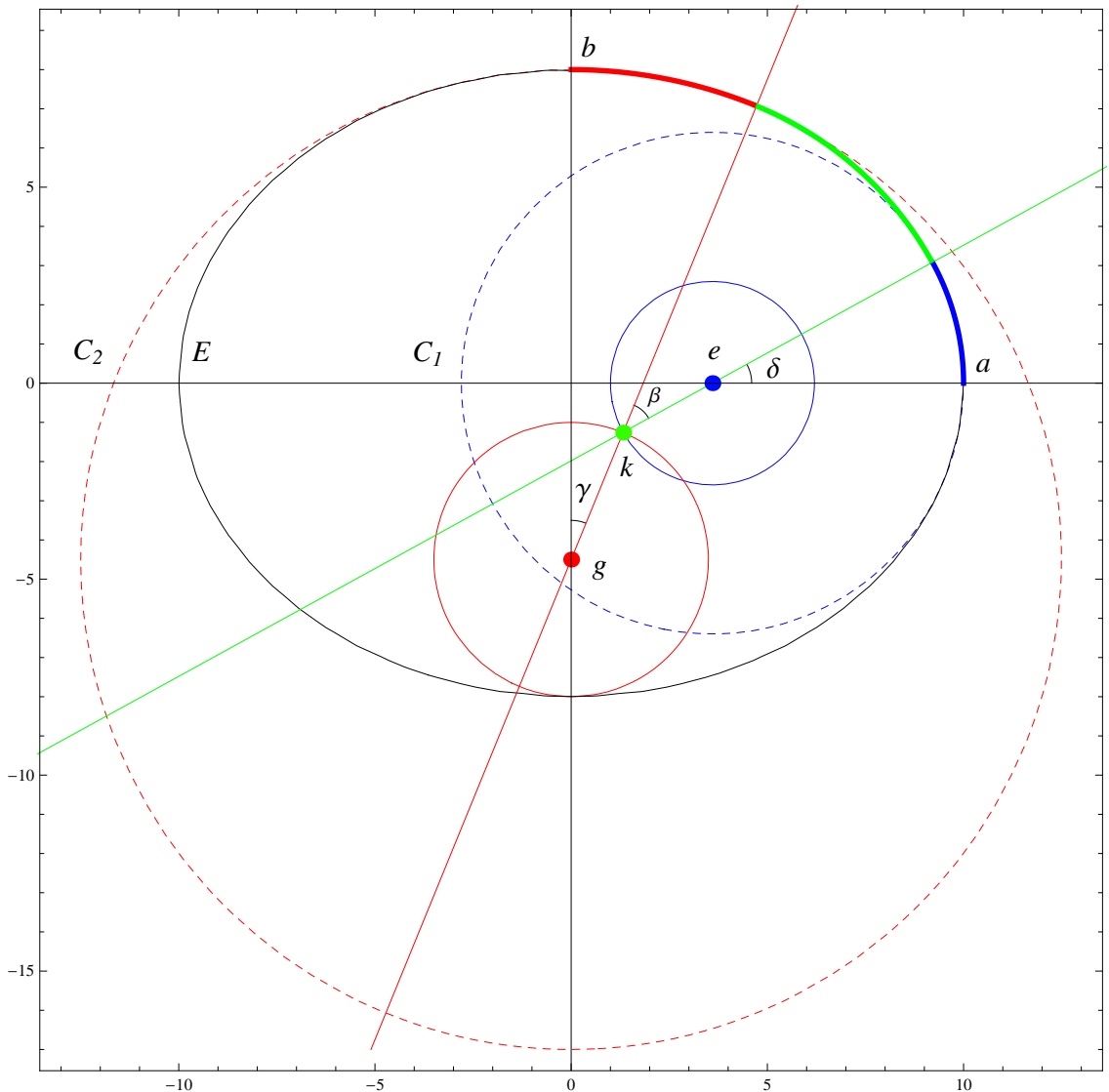}}
\vspace{0.1in}
\caption{Centres et rayons des cercles \textit{mineur}, \textit{interm\'{e}diaire} et \textit{majeur}.}
\label{Fig1}
\end{figure}

On remarque sur la Fig. 1 que la d\'{e}termination du p\'{e}rim\`{e}tre de l'ovale \`{a} huit centres n\'{e}cessite d'expliciter l'expression de la valeur des trois angles aux centres $\gamma$, $\beta$ et $\delta$  en fonction des longueurs du demi-axe majeur $a$ et du demi-axe mineur $b$ de l'ellipse. Pour y parvenir nous proposons d'employer une m\'{e}thode dont les diff\'{e}rentes \'{e}tapes sont d\'{e}crites ci-dessous (voir Fig. 2).

\newpage

\begin{enumerate}
\item[1.] Calcul de l'angle $\alpha'$ dans le triangle ($gek$) avec la ``loi des cosinus''.
\item[2.] Calcul de l'angle $\gamma = \theta - \alpha'$ avec $\theta = \hat{g}$ dans le triangle rectangle ($geo$).
\item[3.] Calcul de l'angle $\alpha$ dans le triangle ($gek$) avec la ``loi des sinus''.
\item[4.] Calcul de l'angle $\beta$ dans le triangle ($gek$) avec la ``loi des sinus''.
\item[5.] Calcul de l'angle $\delta = \theta' - \alpha$ avec $\theta' = \hat{e}$ dans le triangle rectangle ($geo$).
\end{enumerate}

\begin{figure}[htbp]
\centerline{\includegraphics[width=13cm,height=13cm]{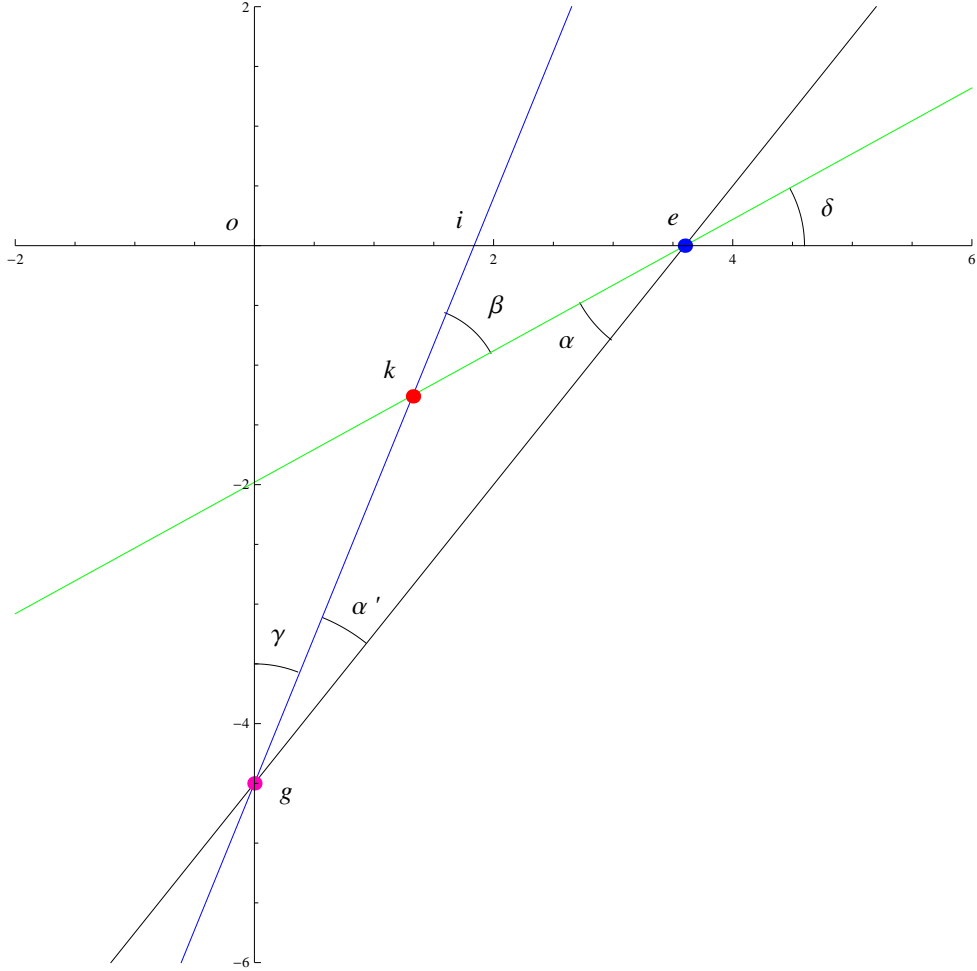}}
\vspace{0.1in}
\caption{Triangles ($geo$) et ($gek$) pour la d\'{e}termination des trois angles aux centres.}
\label{Fig2}
\end{figure}

Au pr\'{e}alable, en appliquant le th\'{e}or\`{e}me de Pythagore au triangle rectangle ($geo$), il est facile de d\'{e}duire l'expression de la longueur de l'hypot\'{e}nuse $|ge|$. En effet, des coordonn\'{e}es des points $g$ et $e$, on d\'{e}duit les distances $|go|$ et $|oe|$. On obtient : $|ge|^2 = |go|^2 + |oe|^2$. Il vient alors :

\begin{equation}
|ge|^2 =\frac{\left(a^2 + b^2\right)\left(a^2 - b^2\right)^2}{a^2 b^2}
\label{eq5}
\end{equation}

\subsection{Calcul de l'angle $\alpha'$ dans le triangle ($gek$) avec la ``loi des cosinus''}

Consid\'{e}rons le triangle ($gek$) (voir Fig. 2) et appliquons la ``loi des cosinus''. On a :

\[
|ke|^2 = |ge|^2 + |gk|^2 - 2|ge|_{}|gk|_{}\cos{\alpha'}
\]

On en d\'{e}duit que:

\begin{equation}
\cos{\left(\alpha'\right)} = \frac{|ge|^2 + |gk|^2 - |ke|^2}{2|ge|_{}|gk|_{}}
\label{eq6}
\end{equation}

Les valeurs de $|gk|$, $|ke|$ et $|ge|$ \'{e}tant donn\'{e}es par les expressions (\ref{eq3}-\ref{eq5}), on obtient finalement :

\begin{equation}
\cos{\left(\alpha'\right)} = \frac{2a^2 + ab + b^2}{\left( b + 2a\right) \sqrt{a^2+b^2}}
\label{eq7}
\end{equation}

De l'identit\'{e} $\cos^2(x)+ \sin^2(x)=1$, on d\'{e}duit la valeur de $\sin(\alpha')$:

\begin{equation}
\sin{\left(\alpha'\right)} = \frac{b}{2a+b}\sqrt{\frac{2ab}{a^2+b^2}}
\label{eq8}
\end{equation}

\subsection{Calcul de l'angle $\gamma = \theta - \alpha'$ avec $\theta = \hat{g}$ dans le triangle rectangle ($geo$)}

Consid\'{e}rons le triangle rectangle ($geo$) (voir Fig. 2) et utilisons les relations m\'{e}triques. On a :

\begin{equation}
\sin{\left(\theta\right)} = \frac{|oe|}{|ge|} = \frac{\left( a^2 - b^2 \right)/a}{\left( a^2 - b^2 \right)\sqrt{a^2+b^2}/ab} = \frac{b}{\sqrt{a^2+b^2}}
\label{eq9}
\end{equation}

On en d\'{e}duit que

\begin{equation}
\sin{\left(\theta\right)} = \frac{b}{\sqrt{a^2+b^2}} \quad \mbox{and} \quad \cos{\left(\theta\right)} = \frac{a}{\sqrt{a^2+b^2}}
\label{eq10}
\end{equation}

De $\gamma  = \theta - \alpha'$, il vient $\sin(\gamma) = \sin(\theta - \alpha') = \sin(\theta)\cos(\alpha') - \cos(\theta)\sin(\alpha')$. En rempla\c{c}ant  $\cos(\alpha')$, $\sin(\theta)$, $\sin(\alpha')$ et $\cos(\theta)$ par les expressions (\ref{eq7}-\ref{eq10}), on obtient :

\begin{equation}
\boxed{\sin{\left(\gamma\right)} = \frac{b}{2a+b} \frac{2a + b + \sqrt{2ab}}{a + b + \sqrt{2ab}}}
\label{eq11}
\end{equation}

\subsection{Calcul de l'angle $\alpha$ dans le triangle ($gek$) avec la ``loi des sinus''}

Consid\'{e}rons de nouveau le triangle ($gek$) (voir Fig. 2) et appliquons la ``loi des sinus''. On a :

\begin{equation}
\frac{|ek|}{\sin{\left(\alpha'\right)}} = \frac{|gk|}{\sin{\left(\alpha\right)}} = \frac{|ge|}{\sin{\left(\pi - \beta\right)}}
\label{eq12}
\end{equation}

Les deux premi\`{e}res parties de cette \'{e}galit\'{e} conduisent \`{a} $\sin{\alpha} = \dfrac{|gk|}{|ek|} \sin(\alpha')$. En rempla\c{c}ant $|gk|$, $|ek|$ et $\sin(\alpha')$ par les expressions (\ref{eq3}, \ref{eq4} \& \ref{eq8}), on obtient :

\begin{equation}
\sin{\left(\alpha\right)} = \frac{a}{a+2b}\sqrt{\frac{2ab}{a^2+b^2}} \quad \mbox{ et } \quad \cos{\left(\alpha\right)} = \frac{a^2+ ab + 2b^2}{\left( a + 2b \right)\sqrt{a^2+b^2}}
\label{eq13}
\end{equation}

\subsection{Calcul de l'angle $\beta$ dans le triangle ($gek$) avec la ``loi des sinus''}

Consider maintenant la premi\`{e}re et derni\`{e}re  partie de l'\'{e}galit\'{e} (\ref{eq12}), on a: $\sin{(\pi - \beta)} = \dfrac{|ge|}{|ek|}\sin(\alpha')$. Alors, en rempla\c{c}ant $|ek|$, $|ge|$ et $\sin(\alpha')$ par les expressions (\ref{eq4}, \ref{eq5} \& \ref{eq8}), et, en tenant compte du fait que : $\sin{(\pi - \beta)} = \sin{(\beta)}$, on obtient :

\begin{equation}
\boxed{\sin{\left(\beta\right)} = \frac{2\left(a + b\right)\sqrt{2ab}}{\left(a + 2b\right)\left(2a + b\right)}}
\label{eq14}
\end{equation}

\subsection{Calcul de l'angle $\delta = \theta' - \alpha$ avec $\theta' = \hat{e}$ dans le triangle rectangle ($geo$)}

Consid\'{e}rons le triangle rectangle ($geo$) (voir Fig. 2), la somme de ces angles est : $\theta + \theta' + \pi/2 = \pi$. Il s'ensuit que $\theta' = \pi/2 - \theta$. Ce qui impliquec que : $\cos{(\theta')} = \cos{(\pi/2 - \theta)} = \sin{(\theta)}$ et $\sin{(\theta')} = \sin{(\pi/2 - \theta)} = \cos{(\theta)}$. Ainsi, \`{a} partir des Eqs. (\ref{eq10}), on obtient :

\begin{equation}
\sin{\left(\theta'\right)} = \cos{\left(\theta\right)} = \frac{a}{\sqrt{a^2+b^2}} \quad \mbox{ et } \quad \cos{\left(\theta'\right)} = \sin{\left(\theta\right)} = \frac{b}{\sqrt{a^2+b^2}}
\label{eq15}
\end{equation}

De $\delta  = \theta' - \alpha$, il vient $\sin(\delta) = \sin(\theta' - \alpha) = \sin(\theta')\cos(\alpha) - \cos(\theta')\sin(\alpha)$. Alors, en rempla\c{c}ant $\cos(\alpha)$, $\sin(\alpha)$, $\sin(\theta')$ et $\cos(\theta')$ par les expressions (\ref{eq13}-\ref{eq15}), on obtient :

\begin{equation}
\boxed{\sin{\left(\delta\right)} = \frac{a}{a + 2b} \frac{a + 2b + \sqrt{2ab}}{a + b + \sqrt{2ab}}}
\label{eq16}
\end{equation}

\textbf{Remarque.}

Il est important de remarquer que les expressions (\ref{eq11}, \ref{eq14} \& \ref{eq16}) pr\'{e}sentent une parfaite sym\'{e}trie. En effet, en rempla\c{c}ant $a$ par $b$ dans l'expression (\ref{eq11}) de $\sin{\left( \gamma \right)}$ on retrouve l'expression (\ref{eq16}) de $\sin{\left( \delta \right)}$ et \textit{vice versa}. De plus, l'expression (\ref{eq14}) de $\sin{\left( \beta \right)}$ est invariante si l'on remplace $a$ par $b$ et \textit{vice versa}. Les principaux r\'{e}sultats obtenus ci-dessus sont r\'{e}sum\'{e}s dans le tableau 1 ci-dessous :

\begin{table}[th!]
\centering
\begin{tabular}{|c|c|c|}
\hline
Cercles & Angles aux centres & Rayons \\
\hline
\hspace{1cm} Majeur \hspace{1cm}  & \hspace{1cm} $\sin{\left(\gamma\right)} = \dfrac{b}{2a+b} \dfrac{2a + b + \sqrt{2ab}}{a + b + \sqrt{2ab}}$ \hspace{1cm} & \hspace{1cm} $R = \dfrac{a^2}{b}$ \hspace{1cm} \\
\hline
\hspace{1cm} Interm\'{e}diaire \hspace{1cm} & \hspace{1cm} $\sin{\left(\beta\right)} = \dfrac{2\left(a + b\right)\sqrt{2ab}}{\left(a + 2b\right)\left(2a + b\right)}$ \hspace{1cm} & \hspace{1cm} $p = \dfrac{a+b}{2}$ \hspace{1cm} \\
\hline
\hspace{1cm} Mineur \hspace{1cm}  & \hspace{1cm} $\sin{\left(\delta\right)} = \dfrac{a}{a + 2b} \dfrac{a + 2b + \sqrt{2ab}}{a + b + \sqrt{2ab}}$ \hspace{1cm} & \hspace{1cm} $r = \dfrac{b^2}{a}$ \hspace{1cm} \\
\hline
\end{tabular}
\vspace{0.25cm}
\caption{Angles aux centres et rayons des cercles \textit{majeur}, \textit{interm\'{e}diaire} et \textit{mineur}.}
\label{tab1}
\end{table}

On peut ainsi construire comme approximation du quart de l'ellipse ($E$) dans le plan $x \geqslant 0$, $y \geqslant 0$, un quart d'ovale form\'{e} de trois arcs de cercles : un arc du cercle \textit{majeur} osculateur \`{a} l'ellipse au point de coordonn\'{e}es ($0,b$), de centre $g$, de rayon $R$ et d'angle au centre $\gamma$ (en rouge sur la Fig. 1); un arc du cercle \textit{mineur} osculateur \`{a} l'ellipse au point de coordonn\'{e}es ($a,0$), de centre $e$, de rayon $r$ et d'angle au centre $\delta$ (en bleu sur la Fig. 1) et un arc du cercle \textit{interm\'{e}diaire}, tangent aux deux autres, de centre $k$, de rayon $p$ et d'angle au centre $\beta$ (en vert sur la Fig. 1).

\newpage

\section{D\'{e}termination du p\'{e}rim\`{e}tre de l'ovale \`{a} huit centres}

En utilisant les expressions des angles et des rayons pr\'{e}sent\'{e}es dans le Tab. 1, le calcul du p\'{e}rim\`{e}tre de l'ovale \`{a} huit centres est donc \'{e}gal \`{a} :

\begin{equation}
\boxed{ O\left( a, b \right) = 4\left[ \arcsin{\left(\gamma\right)}\dfrac{a^2}{b} + \arcsin{\left(\beta\right)}\dfrac{a+b}{2} + \arcsin{\left(\delta\right)}\dfrac{b^2}{a} \right]}
\label{eq17}
\end{equation}

Dans le cas de l'amphith\'{e}\^{a}tre Flavien, plus connu sous le nom de Colis\'{e}e, les dimensions de ce monument fournies par Golvin \cite{Golvin}, Trevisan \cite{Trevisan} et Wilson \cite{Wilson} ainsi que celles relev\'{e}es sur site sont \'{e}gales \`{a} 188 m pour la plus grande longueur et 156 m, pour sa plus grande largeur. En supposant que ce monument a la forme d'une ellipse de demi grand axe $a = 94$  et de demi petit axe $b=78$, la valeur de l'excentricit\'{e} est $\varepsilon = 0.55$. En utilisant l'expression (\ref{eq2}), on obtient pour la valeur du p\'{e}rim\`{e}tre du Colis\'{e}e $L(94,78) = 541.524m$. En utilisant l'expression (\ref{eq17}) du p\'{e}rim\`{e}tre de l'ovale \`{a} huit centres, on obtient alors $O(94,78)=541.523$. L'erreur relative entre ces deux expressions est dans le cas du Colis\'{e}e de l'ordre de $1.85 \times 10^{-4} \%$. Des tests num\'{e}riques ont \'{e}t\'{e} r\'{e}alis\'{e}s afin de calculer les expressions (\ref{eq2}) et (\ref{eq17}) pour diff\'{e}rentes valeurs de $a$ et $b$ variant entre 1 et 10. Ces tests qui pourraient facilement \^{e}tre \'{e}tendus \`{a} des intervalles plus importants, ont montr\'{e} que l'erreur relative entre les deux expressions  (\ref{eq2}) et (\ref{eq17}) reste inf\'{e}rieure \`{a} $0.029\%$.\\

\textbf{Remarque.}

En posant $a = b = \rho$, l'expression (\ref{eq17}) permet de retrouver le p\'{e}rim\`{e}tre du cercle de rayon $\rho$. L'expression (\ref{eq17}) s'\'{e}crit alors :

\begin{equation}
O\left( a, b \right) = 4 \rho \left[ \arcsin{\left(\gamma\right)} + \arcsin{\left(\beta\right)} + \arcsin{\left(\delta\right)} \right]
\label{eq18}
\end{equation}

Dans ce cas, les trois sinus (\ref{eq11}, \ref{eq14} and \ref{eq16}) deviennent \'{e}gaux \`{a} $\sin{\left( \gamma \right)} = \dfrac{4-\sqrt{2}}{6}$, $\sin{\left( \beta \right)} = \dfrac{4\sqrt{2}}{9}$ et $\sin{\left( \delta \right)} = \dfrac{4-\sqrt{2}}{6}$. En rempla\c{c}ant ces sinus dans l'expression (\ref{eq18}) on a :

\begin{equation}
O\left( a, b \right) = 4 \rho \left[ 2\arcsin{(\frac{4-\sqrt{2}}{6})} + \arcsin{(\frac{4\sqrt{2}}{9})} \right]
\label{eq19}
\end{equation}

Alors, en utilisant les formules d'addition des fonctions Arcsinus fournies par Bronstein et Semendiaev \cite[p. 281-282]{Bronstein}, on obtient :

\[
2\arcsin{(\frac{4-\sqrt{2}}{6})} = 2\arcsin{(\frac{7}{9})} \quad \mbox{ et } \quad  \arcsin{(\frac{7}{9})} + \arcsin{(\frac{4\sqrt{2}}{9})} = \pi - \arcsin{(1)} = \frac{\pi}{2}
\]

Ainsi, en rempla\c{c}ant dans l'expression (\ref{eq19}), on d\'{e}montre que lorsque $a = b = \rho$, $O\left( a, b \right) = 2 \pi \rho$.

\section{Conclusion}

On peut donc en conclure d'une part que l'ovale \`{a} huit centres co\"{i}ncide presque parfaitement avec l'ellipse construite sur les m\^{e}mes axes et peut \^{e}tre consid\'{e}r\'{e} comme une repr\'{e}sentation de celle-ci. Et, d'autre part, que le p\'{e}rim\`{e}tre de l'ovale \`{a} huit centres fournit une excellente approximation de celui de l'ellipse construite sur les m\^{e}mes axes.

% Les remerciements sont dans une section, sans num\`{e}rotation

%\section*{Remerciements}
% Remerciements - texte ici

\end{document}